\documentclass[aip,
amsmath,amssymb,
preprint,%
numerical,
]{revtex4-1}
\usepackage{graphicx}
\usepackage{dcolumn}
\usepackage{bm}
\usepackage{hyperref}
\usepackage[utf8]{inputenc}
\usepackage[T1]{fontenc}
\usepackage{mathptmx}
\usepackage{etoolbox}
\usepackage{multirow}
\makeatletter
\def\@email#1#2{%
	\endgroup
	\patchcmd{\titleblock@produce}
	{\frontmatter@RRAPformat}
	{\frontmatter@RRAPformat{\produce@RRAP{*#1\href{mailto:#2}{#2}}}\frontmatter@RRAPformat}
	{}{}
}%
\makeatother
\begin{document}

	\title{Effect of Rotation in an Isotropic Scattering Algal Suspension with Oblique Collimated Irradiation}
	\author{S. K. Rajput}
	\altaffiliation[Corresponding author email: ]{shubh.iiitj@gmail.com}
	\affiliation{ 
		Department of Mathematics, PDPM Indian Institute of Information Technology Design and Manufacturing,
		Jabalpur 482005, India.
	}%
	
	

	\begin{abstract}
		The linear stability of a suspension of isotropic scattering phototactic algae is investigated numerically with particular emphasis on the effects of Taylor number in the rotating medium. The suspension is illuminated by the oblique collimated irradiation. The solutions show a transition of the most unstable mode from stationary to an overstable state or vice versa for certain parameters at the variation in the Taylor number. Oscillatory instabilities are also observed at the three-quarter height of the suspension for some parameters. 
		
	\end{abstract}
	
	\maketitle

	\section{INTRODUCTION}
	The $bioconvection$ is a collective motion of a fluid that holds self-propelled motile microorganisms. These microorganisms, which are generally denser than the surrounding medium (such as water), tend to move upward on average within it. Primarily, this behavior is displayed by algae and bacteria. The arrangement of patterns in bioconvection vanishes when these microorganisms cease their motion. Nevertheless, it's noteworthy that pattern creation does not exclusively depend on the combination of upward motion and higher density. Instances exist where neither upward motion nor higher density contributes to pattern formation. The orientation of these swimming microorganisms is influenced by various environmental cues, referred to as $taxes$. $Gravitaxis$, $chemotaxis$, $phototaxis$, and $gyrotaxis$ are notable forms of $taxes$. $Gravitaxis$ refers to the response to gravitational force, $chemotaxis$ corresponds to the reaction to chemicals, $gyrotaxis$ arises from the balance between gravity-induced torque and local shear flow, and $phototaxis$ defines movement either toward (positive phototaxis) or away from (negative phototaxis) a light source. This discussion solely focuses on the impact of phototaxis.
	
	Experimental observations underscore that the arrangement of patterns in bioconvection can be substantially influenced by diverse degrees of illuminative intensity (e.g., diffuse irradiation)\cite{1wager1911,2kitsunezaki2007}. Intense light can disrupt stable patterns or hinder their formation in a suspension of motile microorganisms within a well-stirred culture. Factors such as light intensity can induce alterations in the size, shape, structure, and symmetry of these patterns\cite{3kessler1985,4williams2011,5kessler1989}. Variations in bioconvection patterns due to light intensity can be elucidated by the following considerations. Firstly, phototactic algae derive energy from photosynthesis, which consequently modifies their swimming behavior (phototaxis). When the light intensity $I$ surpasses or falls below its critical value $I_c$, cells exhibit movement towards (positive phototaxis) or away from (negative phototaxis) the light source. This leads the algae cells to aggregate in locations offering the optimal light intensity ($I=I_c$). Another plausible reason for pattern alterations may be the absorption and scattering of light. The cells' absorption of light causes a reduction in light intensity along the incident path, while scattering initially diverts light away from the path, resulting in decreased intensity, followed by scattering from other points that lead to intensity enhancement~\cite{7ghorai2010}.
	
	In this particular investigation, we make use of the phototaxis model introduced by Rajput. The setup involves a collection of phototactic microorganisms within a rotating medium, revolving around the z-axis at a consistent angular velocity. These microorganisms are exposed to directed irradiation from above. Given the reliance of many motile algae on photosynthesis for sustenance, they exhibit noticeable phototaxis. To accurately depict their actions, it becomes essential to scrutinize the phototaxis model in the context of a rotating medium.
	
	However, the original Rajput model overlooked the influence of scattering prompted by the microorganisms. Consequently, these microorganisms only received direct light from the source directly overhead, resulting in a purely vertical swimming orientation. The velocity of their movement was established based on the amount of light reaching each cell. In our revised model, we integrate the scattering impact generated by the cells. This signifies that the microorganisms now engage with light through scattering, introducing a more intricate swimming behavior. Consequently, the microorganisms' movement direction comprises not only a vertical component but also a horizontal one. By incorporating scattering, our model offers a more precise portrayal of how microorganisms react to light.
	
	Considering an oblique collimated irradiation for the illuminated suspension, equilibrium arises when the fluid's velocity reaches zero. The vertical motion driven by positive and negative phototaxis is counterbalanced by diffusion arising from the random movement of cells. Consequently, a concentrated sublayer emerges, with its position determined by the critical intensity $I_c$. If the intensity across the suspension surpasses (falls below) the critical intensity, a sublayer emerges above (below) the suspension. If the critical intensity $I_c$ falls between the maximum and minimum values of light intensity, a sublayer forms between the upper and lower boundaries of the suspension. The upper sublayer represents a gravitationally stable zone, while the lower sublayer is unstable gravitationally. If the fluid layer becomes unstable, it transitions from the lower unstable region to infiltrate the upper stable region, illustrating a case of penetrative bioconvection~\cite{9straughan1993,10ghorai2005,11panda2016}.
	
	\begin{figure}[!h]
		\centering
		\includegraphics[width=14cm]{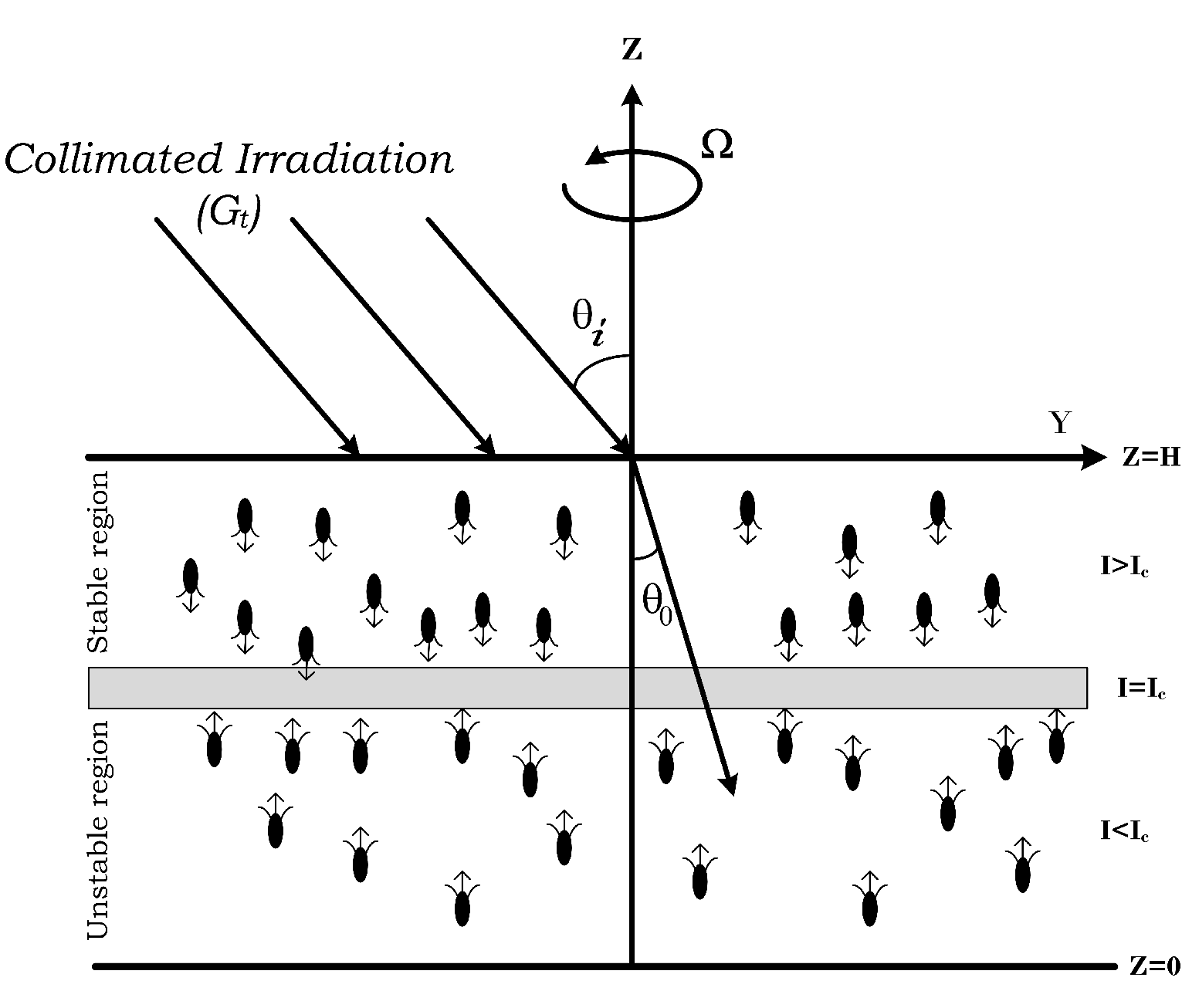}
		\caption{\footnotesize{Accumulation of algae cells at $I=I_c$ in the rotating medium.}}
		\label{fig1}
	\end{figure}
	
	Vincent and Hill~\cite{12vincent1996} endeavoured to simulate phototactic bioconvection for algae suspension in a shallow layer, focusing on cases where there is no orientation bias due to gravity. They developed a simple top-down photoresponse mechanism based on light intensity. Their emphasis rested on the cells' negative buoyancy, suggesting that other overall stress sources were neglected. Utilizing the well-known Beer-Lambert law to model light intensity, they performed a linear stability analysis, uncovering the stationary and oscillatory nature of disturbances. Ghorai and Hill~\cite{10ghorai2005} delved into the quantification of two-dimensional phototactic bioconvection within a confined suspension, bounded by a rigid bottom and stress-free upper and side walls. They addressed significant inaccuracies in Vincent and Hill's equilibrium solution. Building upon this, Ghorai $et$ $al$.~\cite{7ghorai2010} explored the impact of light scattering, revealing bimodal equilibrium-state configurations for suspensions characterized by substantial scattering. Further research ventured into different aspects. Ghorai and Panda~\cite{13ghorai2013} utilized linear theory to scrutinize bioconvection initiation in an anisotropic scattering solution of phototactic algae, emphasizing forward scattering effects. Panda and Ghorai~\cite{14panda2013} simulated two-dimensional phototactic bioconvection in an isotropic scattering suspension through linear simulations, showing qualitative distinctions in bioconvective patterns compared to Ghorai and Hill~\cite{10ghorai2005} at higher critical wavelengths. Panda and Singh~\cite{11panda2016} numerically investigated linear stability using Vincent and Hill's continuum model, observing a considerable stabilizing effect on suspension due to rigid side walls. However, these prior studies did not account for the effects of diffuse irradiation. Addressing this gap, Panda $et$ $al$.\cite{15panda2016} introduced a model examining the impact of diffuse radiation on isotropic scattering algal suspensions, discovering strong stabilization effects. They also noted that diffuse irradiation had a greater influence on the critical state compared to collimated irradiation alone. Panda\cite{8panda2020} extended the investigation to include forward anisotropic scattering's impact on phototactic bioconvection under both diffuse and collimated irradiation. Additionally, Kumar~\cite{17kumar2022} explored the effect of oblique collimated irradiation on isotropic scattering algal suspensions, finding various solution types for specific parameter ranges. Recently, research has examined the combined impact of oblique and diffuse irradiation on bioconvection onset~\cite{41rajput2023}, and Rajput investigated the stability of rotating algal suspensions under oblique irradiation, focusing on absorption effects without considering scattering. This opens the door to further investigation. In light of this, we propose a model for isotropic scattering algal suspensions subjected to rotation and illuminated by oblique collimated irradiation.
	
	The structure of this article is organized as follows: We commence by presenting a mathematical formulation of the phototaxis model and determining the basic equilibrium state. Subsequently, we analyze the linear stability of this equilibrium state, solving it numerically to derive neutral curves. Finally, we discuss the obtained results using these neutral curves.
	
	
	\section{MATHEMATICAL FORMULATION}
	
	Consider the phototactic algal suspension within a layer of finite depth $H$ and infinite lateral extent. In this model, the suspension is illuminated from the top by oblique collimated irradiation with the magnitude $G_t$.
	
	
	\subsection{\label{sec:level3}PHOTOTAXIS WITH ABSORPTION AND SCATTERING}
	
	In the absorbing and scattering medium, the light intensity profile is calculated by radiative transfer equation (RTE)
	
	\begin{equation}\label{1}
		\frac{dG(\boldsymbol{x},\boldsymbol{s})}{ds}+(a+\sigma_s)G(\boldsymbol{x},\boldsymbol{s})=\frac{\sigma_s}{4\pi}\int_{0}^{4\pi}G(\boldsymbol{x},\boldsymbol{s'})\Phi(\boldsymbol{s},\boldsymbol{s'})d\Omega',
	\end{equation}
	
	where $G(\boldsymbol{x},\boldsymbol{s})$ is the light inetensity  at $\boldsymbol{x}$ in $\boldsymbol{s}$ direction. $\alpha and \sigma_s$ are the absorption and scattering coefficients, respectively, and $\Phi(\boldsymbol{s},\boldsymbol{s'})$ is the scattering phase function. Now, assume that the algae cells are scattered light isotropically, i.e., $\Phi(\boldsymbol{s},\boldsymbol{s'})=1$~\cite{15panda2016}.
	
	\begin{figure}[!h]
		\centering
		\includegraphics[width=14cm]{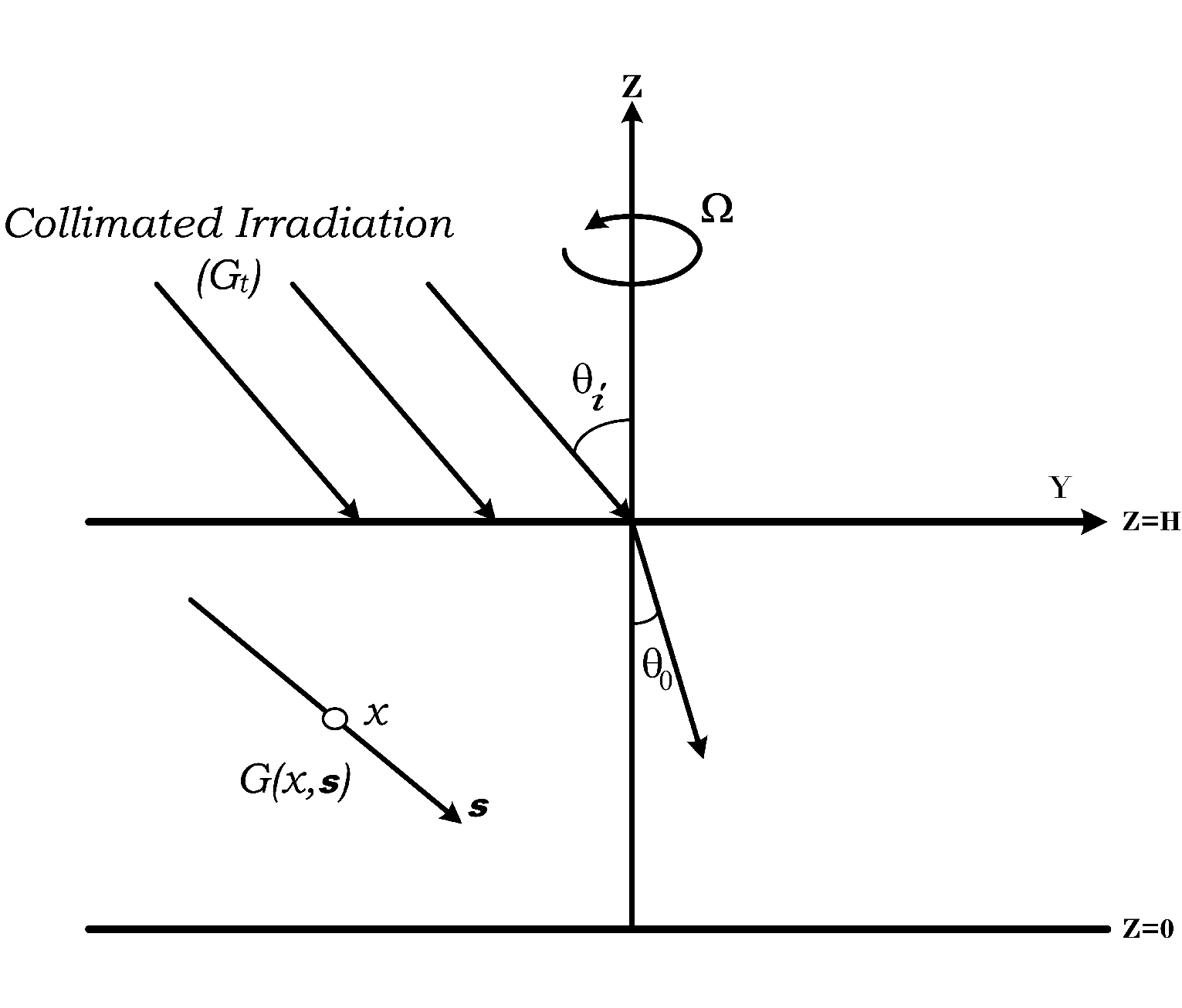}
		\caption{\footnotesize{Axial representation of the problem.}}
		\label{fig2}
	\end{figure}

	We assume  that the suspension's boundaries are non-reflective, then the intensity at the top of the suspension is thus given by
	\begin{equation*}
		G(\boldsymbol{x}_b,\boldsymbol{s})=G_t\delta(\boldsymbol{s}-\boldsymbol{s_0}), 
	\end{equation*}
	where $\boldsymbol{x}_b=(x,y,H)$ is the location at the top boundary surface and $\boldsymbol{s_0}=\sin(\pi-\theta_0)[\cos\phi+\sin\phi]+\cos(\pi-\theta_0)$. If we take $\theta_i=0$, this model will be similar to the previous model proposed by Rajput. Consider, $a=\alpha n(\boldsymbol{x})$ and $\sigma_s=\beta n(\boldsymbol{ x})$, then the RTE becomes
	\begin{equation}\label{2}
		\frac{dG(\boldsymbol{x},\boldsymbol{s})}{ds}+(\alpha+\beta)nG(\boldsymbol{x},\boldsymbol{s})=\frac{\beta n}{4\pi}\int_{0}^{4\pi}G(\boldsymbol{x},\boldsymbol{s'})d\Omega'.
	\end{equation}
	In the medium, the total intensity at a fixed point $\boldsymbol{x}$ is 
	\begin{equation*}
		I(\boldsymbol{x})=\int_0^{4\pi}G(\boldsymbol{x},\boldsymbol{s})d\Omega,
	\end{equation*}
	and the radiative heat flux is given by
	\begin{equation}\label{3}
		\boldsymbol{q}(\boldsymbol{x})=\int_0^{4\pi}G(\boldsymbol{x},\boldsymbol{s})\boldsymbol{s}d\Omega.
	\end{equation}
	Let $\boldsymbol{p}$ be the unit vector in the cell swimming direction and $<\boldsymbol{p}>$ is the mean swimming direction.
	For many species of microorganisms, the swimming speed is free from the light condition, position, time, and direction~\cite{18hill1997}. Here, we assumed that the cells and the fluid flow at the same speed, then the average cell
	Swimming velocity is defined as
	\begin{equation*}
		\boldsymbol{W}_c=W_c<\boldsymbol{p}>,
	\end{equation*}
	where $W_c$ is the average cell swimming speed and the cell's mean swimming direction $<\boldsymbol{p}>$ is calculated by
	\begin{equation}\label{4}
		<\boldsymbol{p}>=-T(I)\frac{\boldsymbol{q}}{|\boldsymbol{q}|},
	\end{equation}
	where $\varpi$ is a non-negative constant and $T(G), $ is the taxis response function (taxis function), which shows the response of algae cells to light and has the mathematical form such that 
	\begin{equation*}
		T(I)=\left\{\begin{array}{ll}\geq 0, & \mbox{ } I(\boldsymbol{x})\leq I_{c}, \\
			< 0, & \mbox{ }I(\boldsymbol{x})>I_{c}.  \end{array}\right. 
	\end{equation*}
	
	The mean swimming direction becomes zero at the critical light intensity ( $G = G_c)$.
	Generally, the exact functional form of taxis function depends on the species of the microorganisms~\cite{12vincent1996}.
	
	
	\subsection{GOVERNING EQUATIONS}
	
	Suppose we're studying a group of similar algal cells that are spread out and can be represented as a continuous distribution. This algal suspension is very dilute $(nv<<1)$. Each cell has a certain volume $v$ and is a bit denser than the fluid it's in $(\Delta\rho<<\rho)$. The fluid is moving on average at velocity $\boldsymbol{u}$, and the $n$ cells are moving within it. The whole suspension is rotating around a z-axis vertically upward with a constant angular velocity $\Omega$. We're assuming that the suspension is incompressible. Then the governing equations
	
	\begin{equation}\label{5}
		\boldsymbol{\nabla}\cdot \boldsymbol{u}=0.
	\end{equation}
	
	\begin{equation}\label{6}
		\rho\frac{D \boldsymbol{u}}{D t}+2\rho(\boldsymbol{\Omega}\times \boldsymbol{u})=-\boldsymbol{\nabla} P_e+\mu\nabla^2\boldsymbol{u}-nvg\Delta\rho\hat{\boldsymbol{z}},
	\end{equation}
	
	where $D/Dt=\partial/\partial t+\boldsymbol{u}\cdot\boldsymbol{\nabla}$ is the material or total derivative, $P_e$ is the excess pressure above hydrostatic, $\mu$ is the viscosity of the suspension which is assumed to be that of fluid.
	
	\begin{equation}\label{7}
		\frac{\partial n}{\partial t}=-\boldsymbol{\nabla}\cdot \boldsymbol{F},
	\end{equation}
	where $\boldsymbol{F}$ is the total cell flux which is given by
	\begin{equation}\label{8}
		\boldsymbol{F}=nu+nW_c<\boldsymbol{p}>-\boldsymbol{D}\boldsymbol{\nabla} n.
	\end{equation}
	
	Here, the diffusivity tensor $\boldsymbol{D}$ is assumed to be isotropic and constant as $\boldsymbol{D} = DI$. $\boldsymbol{F}$ has two key assumptions. First, cells are assumed to be purely phototactic. Second, a constant value for the diffusion tensor is assumed. Making these assumptions may remove the Fokker-Planck equation from the governing equations.
	
	Consider, that the lower boundary is rigid and the top boundary is stress-free. Then the boundary conditions
	\begin{equation}\label{9}
		\boldsymbol{u}\cdot\hat{\boldsymbol{z}}=\boldsymbol{u}\times\hat{\boldsymbol{z}}=\boldsymbol{F}\cdot\hat{\boldsymbol{z}}=0\qquad on\quad z=0,
	\end{equation}
	\begin{equation}\label{10}
		\boldsymbol{u}\cdot\hat{\boldsymbol{z}}=\frac{\partial^2}{\partial z^2}(\boldsymbol{u}\cdot\hat{\boldsymbol{z}})=\boldsymbol{F}\cdot\hat{\boldsymbol{z}}=0\qquad on\quad z=H.
	\end{equation}

	Boundary conditions for the light intensities
	\begin{subequations}
		\begin{equation}\label{12a}
			G(x, y, z = 0, \theta, \phi)=G_t\delta(\boldsymbol{s}-\boldsymbol{s_0}),
		\end{equation}
		\begin{equation}\label{12b}
			G(x, y, z = 0, \theta, \phi) =0,\quad (0\leq\theta\leq\pi/2).
		\end{equation}
	\end{subequations}

	The governing equations are made dimensionless
	by choosing $H$ as length-scale, $H^2/D$ as time scale, $D/H$ as velocity scale, $\mu D/H^2$ as pressure scale, and $\Bar{n}$ as concentration scale, which are given below 
	\begin{equation}\label{12}
		\boldsymbol{\nabla}\cdot\boldsymbol{u}=0,
	\end{equation}
	\begin{equation}\label{13}
		Sc^{-1}\left(\frac{\partial \boldsymbol{u}}{\partial t}+(\boldsymbol{u}\cdot\nabla )\boldsymbol{u}\right)+\sqrt{Ta}(\hat{z}\times\boldsymbol{u})=-\nabla P_{e}-Rn\hat{\boldsymbol{z}}+\nabla^{2}\boldsymbol{u},
	\end{equation}
	\begin{equation}\label{14}
		\frac{\partial{n}}{\partial{t}}=-{\boldsymbol{\nabla}}\cdot{\boldsymbol{F}},
	\end{equation}
	where
	\begin{equation}\label{15}
		{\boldsymbol{F}}=n{\boldsymbol{u}}+nV_{c}<{\boldsymbol{p}}>-{\boldsymbol{\nabla}}n.
	\end{equation}

	Here, $S{c}=\nu/{D}$ is the Schmidt number, $V_c=W_cH/D$ is scaled swimming speed, $R=\overline{n}\vartheta g\Delta{\rho}H^{3}/\nu\rho{D}$ is the Rayleigh number, and $Ta=4\Omega^2H^4/\nu^2$ is the Taylor number.
	In dimensionless form, the boundary conditions become
	
	\begin{equation}
		\boldsymbol{u}\cdot\hat{\boldsymbol{z}}=\boldsymbol{u}\times\hat{\boldsymbol{z}}=\boldsymbol{F}\cdot\hat{\boldsymbol{z}}=0\qquad on\quad z=0,
	\end{equation}
	\begin{equation}
		\boldsymbol{u}\cdot\hat{\boldsymbol{z}}=\frac{\partial^2}{\partial z^2}(\boldsymbol{u}\cdot\hat{\boldsymbol{z}})=\boldsymbol{F}\cdot\hat{\boldsymbol{z}}=0\qquad on\quad z=1.
	\end{equation}
	
	The RTE in dimensionless form
	
	\begin{equation}\label{18}
		\xi\frac{dG}{dx}+\eta\frac{dG}{dy}+\nu\frac{dG}{dz}+\kappa nG(\boldsymbol{x},\boldsymbol{s})=\frac{\omega\kappa n}{4\pi}\int_{0}^{4\pi}G(\boldsymbol{x},\boldsymbol{s'})d\Omega',
	\end{equation}
	where $\xi,\eta$ and $\nu$ are the direction cosines, $\kappa=(\alpha+\beta)\Bar{n}H$, $\sigma=\beta\Bar{n}H$, and $\omega=\sigma/\kappa$.
	
	In dimensionless form, the intensity at boundaries becomes,
	\begin{subequations}
		\begin{equation}
			G(x, y, z = 1, \theta, \phi)=G_t\delta(\boldsymbol{s}-\boldsymbol{s_0}),
		\end{equation}
		\begin{equation}
			G(x, y, z = 0, \theta, \phi) =0,\qquad (0\leq\theta\leq\pi/2). 
		\end{equation}
	\end{subequations}
	
	
	\section{THE BASIC (EQUILIBRIUM) STATE SOLUTION}
	
	Equations $(\ref{12})-(\ref{14})$ and $(\ref{18})$ have an equilibrium solution in the form of
	
	\begin{equation}\label{24}
		\boldsymbol{u}=0,~~~n=n_s(z)\quad and\quad  G=G_s(z,\theta).
	\end{equation}
	Therefore, the total intensity $I_s$ and radiative flux $\boldsymbol{q}_s$ at the equilibrium state are given by
	
	\begin{equation*}
		I_s=\int_0^{4\pi}G_s(z,\theta)d\Omega,\quad 
		\boldsymbol{q}_s=\int_0^{4\pi}G_s(z,\theta)\boldsymbol{s}d\Omega,
	\end{equation*}
	and the governing equation for $G_s$ can be written as
	\begin{equation}\label{25}
		\frac{dG_s}{dz}+\frac{\kappa n_sG_s}{\nu}=\frac{\omega\kappa n_s}{4\pi\nu}I_s(z).
	\end{equation}
	
	We decompose the basic state intensity into collimated part $G_s^c$ and diffuse part $G_s^d$ ( which occurs due to scattering) such that $G_s=G_s^c+G_s^d$. The equation governs the collimated part $I_s^c$. 
	
	\begin{equation}\label{26}
		\frac{dG_s^c}{dz}+\frac{\kappa n_sG_s^c}{\nu}=0,
	\end{equation}
	
	subject to the boundary conditions

	\begin{equation}\label{27}
		G_s^c( 1, \theta) =G_t\delta(\boldsymbol{s}-\boldsymbol{ s}_0),\qquad (\pi/2\leq\theta\leq\pi), 
	\end{equation}
	Now $G_s^c$ is given by
	\begin{equation}\label{28}
		G_s^c=G_t\exp\left(\int_z^1\frac{\kappa n_s(z')}{\nu}dz'\right)\delta(\boldsymbol{s}-\boldsymbol{s_0}), 
	\end{equation}
	
	and the diffused part is governed by  
	\begin{equation}\label{29}
		\frac{dG_s^d}{dz}+\frac{\kappa n_sG_s^d}{\nu}=\frac{\omega\kappa n_s}{4\pi\nu}I_s(z),
	\end{equation}
	subject to the boundary conditions
	\begin{subequations}
		\begin{equation}\label{30a}
			G_s^d( 1, \theta) =0,\qquad (\pi/2\leq\theta\leq\pi), 
		\end{equation}
		\begin{equation}\label{30b}
			G_s^d( 0, \theta) =0,\qquad (0\leq\theta\leq\pi/2). 
		\end{equation}
	\end{subequations}

	Now the total intensity, $I_s=I_s^c+I_s^d$ in the equilibrium state can be written as
	\begin{equation}\label{31}
		I_s^c=\int_0^{4\pi}G_s^c(z,\theta)d\Omega=G_t\exp\left(\frac{-\int_z^1\kappa n_s(z')dz'}{\cos\theta_0}\right),
	\end{equation}
	\begin{equation}\label{32}
		I_s^d=\int_0^{\pi}G_s^d(z,\theta)d\Omega.
	\end{equation}
	
	We get the well-known Lambert-Beer law $I_s=I_s^c$ for no scattering. If we define the optical thickness as 
	\begin{equation*}
		\tau=\int_z^1 \kappa n_s(z')dz',
	\end{equation*}
	
	then the total intensity $I_s$ becomes a function of $\tau$ only. Further, the non-dimensional total intensity, $\Lambda(\tau)=I_s(\tau)/G_t$, satisfies the following Fredholm Integral Equation (FIE),
	
	\begin{equation}\label{33}
		\Lambda(\tau) = e^{-\tau/\cos\theta_0}+\frac{\omega}{2}\int_0^\kappa \Lambda(\tau')E_1(|\tau-\tau'|)d\tau',
	\end{equation}
	
	Where $E_n(x)$ is the exponential integral of order n. This FIE has a singularity at $\tau'=\tau$, therefor solve this FIE, the method of subtraction of singularity is utilised
	
	The basic state radiative flux is written as
	
	\begin{equation*}
		\boldsymbol{q_s}=\int_0^{4\pi}\left(G_s^c(z,\theta)+G_s^d(z,\theta)\right)\boldsymbol{s}d\Omega=-G_t(\cos\theta_0)\exp\left(\frac{\int_1^z\kappa n_s(z')dz'}{\cos\theta_0}\right)\hat{\boldsymbol{z}}+\int_0^{4\pi}G_s^d(z,\theta)\boldsymbol{s}d\Omega.
	\end{equation*}
	
	Here, $\boldsymbol{q_s}$ is free from the x and y component. Therefore, $\boldsymbol{q}_s=-q_s\hat{\boldsymbol{z}}$, where $q_s=|\boldsymbol{q_s}|$. The mean swimming direction becomes
	
	\begin{equation*}
		<\boldsymbol{p_s}>=-T(I_s)\frac{\boldsymbol{q_s}}{q_s}=T(I_s)\hat{\boldsymbol{z}},
	\end{equation*}
	
	The basic cell concentration $n_s(z)$ satisfies,
	
	\begin{equation}\label{34}
		\frac{dn_s}{dz}=V_cT(I_s)n_s,
	\end{equation}
	with the relation
	\begin{equation}\label{35}
		\int_0^1n_s(z)dz=1.
	\end{equation}
	Eqs.~(\ref{33})-(\ref{35}) constitute a boundary value problem which is solved numerically by using a
	shooting method.
	
	\section{Linear stability of the problem}
	
	We consider a small perturbation $\epsilon (0<\epsilon<<1)$ to the equilibrium state as
	
	\begin{widetext}
		\begin{align}\label{37}
			\nonumber[\boldsymbol{u},\zeta,n,G,<\boldsymbol{p}>]=[0,\zeta_s,n_s,G_s,<\boldsymbol{p}_s>]+{\epsilon} [\boldsymbol{u}_1,\zeta_1,n_1,G_1,<\boldsymbol{p}_1>]+{O}({\epsilon}^2). 
		\end{align}
	\end{widetext}
	
	The perturbed variables are substituted into the Eqs.~(\ref{12})-(\ref14) and linearizing about the equilibrium state by collecting $O(\epsilon)$ terms, gives
	
	\begin{equation}\label{39}
		\boldsymbol{\nabla}\cdot \boldsymbol{u}_1=0,
	\end{equation}
	
	where  $\boldsymbol{u}_1=(u_1,v_1,w_1)$.
	
	\begin{equation}\label{40}
		Sc^{-1}\left(\frac{\partial \boldsymbol{u_1}}{\partial t}\right)+\sqrt{Ta}(z\times \boldsymbol{u}_1)+\boldsymbol{\nabla} P_{e}+Rn_1\hat{\boldsymbol{z}}=\nabla^{2}\boldsymbol{ u_1},
	\end{equation}
	
	\begin{equation}\label{41}
		\frac{\partial{n_1}}{\partial{t}}+V_c\boldsymbol{\nabla}\cdot(<\boldsymbol{p_s}>n_1+<\boldsymbol{p_1}>n_s)+w_1\frac{dn_s}{dz}=\boldsymbol{\nabla}^2n_1.
	\end{equation}
	
	If $I=I_s+\epsilon I_1+{O}(\epsilon^2)=(I_s^c+\epsilon I_1^c)+(I_s^d+\epsilon I_1^d)+{O}(\epsilon^2)$, then the steady collimated total intensity is perturbed as $G_t\exp\left(-\kappa\int_z^1(n_s(z')+\epsilon n_1+\mathcal{O}(\epsilon^2))dz'\right)$  and after simplification, we get
	
	\begin{equation}\label{42}
		I_1^c=G_t\exp\left(\frac{-\int_z^1 \kappa n_s(z')dz'}{\cos\theta_0}\right)\left(\frac{\int_1^z\kappa n_1 dz'}{\cos\theta_0}\right)
	\end{equation}
	
	Similarly, $I_1^d$ is given by
	
	\begin{equation}\label{43}
		I_1^d=\int_0^{4\pi}G_1^d(\boldsymbol{ x},\boldsymbol{ s})d\Omega.
	\end{equation}
	
	Similarly, for the radiative heat flux $q=q_s+\epsilon q_1+{O}(\epsilon^2)=(q_s^c+\epsilon q_1^c)+(q_s^d+\epsilon q_1^d)+{O}(\epsilon^2)$, we find 
	
	\begin{equation}\label{44}
		\boldsymbol{q}_1^c=-G_t(\cos\theta_0)\exp\left(\frac{-\int_z^1 \kappa n_s(z')dz'}{\cos\theta_0}\right)\left(\frac{\int_1^z\kappa n_1 dz'}{\cos\theta_r}\right)\hat{z}
	\end{equation}
	
	and
	
	\begin{equation}\label{45}
		q_1^d=\int_0^{4\pi}G_1^d(\boldsymbol{ x},\boldsymbol{ s})\boldsymbol{ s}d\Omega.
	\end{equation}
	
	Now the expression
	
	\begin{equation*}
		-T(I_s+\epsilon I_1)\frac{\boldsymbol{q}_s+\epsilon\boldsymbol{q}_1+{O}(\epsilon^2)}{|\boldsymbol{q}_s+\epsilon\boldsymbol{q}_1+{O}(\epsilon^2)|}-T(I_s)\hat{\boldsymbol{z}},
	\end{equation*}
	
	gives the perturbed swimming direction on collecting $O(\epsilon)$ terms
	
	\begin{equation}\label{46}
		<\boldsymbol{p_1}>=I_1\frac{dT(I_s)}{dI}\hat{\boldsymbol{z}}-T(I_s)\frac{\boldsymbol{q_1}^H}{\boldsymbol{q_s}},
	\end{equation}
	
	where $\boldsymbol{q}_1^H=[\boldsymbol{q}_1^x,\boldsymbol{q}_1^y]$ is the horizontal component of the perturbed radiative flux $\boldsymbol{q}_1$.
	Now substituting the value of $<\boldsymbol{p_1}>$ from  Eq.~$(\ref{46})$ into Eq.~$(\ref{41})$ and simplifying, we get
	
	\begin{equation}\label{47}
		\frac{\partial{n_1}}{\partial{t}}+V_c\frac{\partial}{\partial z}\left(T(I_s)n_1+n_sI_1\frac{dT(I_s)}{dI}\right)-V_cn_s\frac{T(I_s)}{q_s}\left(\frac{\partial q_1^x}{\partial x}+\frac{\partial q_1^y}{\partial y}\right)+w_1\frac{dn_s}{dz}=\nabla^2n_1.
	\end{equation}
	
	By elimination of $P_e$ and the horizontal component of $u_1$, Eqs. (26), (27) and Eq. (31) can be reduced to three equations for the perturbed variables namely the vertical component of the velocity $w_1$, the vertical component of the vorticity $\zeta_1 (= \zeta\cdot\hat{\boldsymbol{z}})$ and the concentration $n_1$. These variables can be decomposed into normal modes as
	
	\begin{equation}\label{47}
		[w_1,\zeta_1,n_1]=[W(z),Z(z),N(z)]\exp{(\sigma t+i(lx+my))}.  
	\end{equation}
	
	The governing equation for perturbed intensity $G_1$ can be written as
	
	\begin{equation}\label{49}
		\xi\frac{\partial G_1}{\partial x}+\eta\frac{\partial G_1}{\partial y}+\nu\frac{\partial G_1}{\partial z}+\kappa( n_sG_1+n_1G_s)=\frac{\omega\kappa}{4\pi}(n_sI_1+I_sn_1),
	\end{equation}
	
	subject to the boundary conditions
	
	\begin{subequations}
		\begin{equation}\label{50a}
			G_1(x, y, z = 1, \xi, \eta, \nu) =0,\qquad (\pi/2\leq\theta\leq\pi,0\leq\phi\leq 2\pi), 
		\end{equation}
		\begin{equation}\label{50b}
			G_1(x, y, z = 0,\xi, \eta, \nu) =0,\qquad (0\leq\theta\leq\pi/2,0\leq\phi\leq 2\pi). 
		\end{equation}
	\end{subequations}
	
	Now
	
	\begin{equation*}
		G_1^d=\Psi^d(z,\xi,\eta,\nu)\exp{(\sigma t+i(lx+my))}. 
	\end{equation*}
	
	From Eqs.~(\ref{42}) and (\ref{43}), we get
	
	\begin{equation}\label{51}
		I_1^c=\left[G_t\exp\left(\frac{\int_1^z \kappa n_s(z')dz'}{\cos\theta_0}\right)\left(\frac{\int_1^z\kappa n_1 dz'}{\cos\theta_0}\right)\right]\exp{(\sigma t+i(lx+my))}=\mathcal{I}^c(z)\exp{(\sigma t+i(lx+my))},
	\end{equation}
	
	and
	
	\begin{equation}\label{52}
		I_1^d=\mathcal{I}^d(z)\exp{(\sigma t+i(lx+my))}= \left(\int_0^{4\pi}\Psi^d(z,\xi,\eta,\nu)d\Omega\right)\exp{(\sigma t+i(lx+my))},
	\end{equation}
	
	where $\mathcal{I}(z)=\mathcal{I}^c(z)+\mathcal{I}^d(z)$ is the perturbed total intensity.
	
	Now $\Psi^d$ satisfies
	
	\begin{equation}\label{53}
		\frac{d\Psi^d}{dz}+\frac{(i(l\xi+m\eta)+\kappa n_s)}{\nu}\Psi^d=\frac{\omega\kappa}{4\pi\nu}(n_s\mathcal{I}+I_s\Theta)-\frac{\kappa}{\nu}G_s\Theta,
	\end{equation}
	
	subject to the boundary conditions
	
	\begin{subequations}
		\begin{equation}\label{54a}
			\Psi^d( 1, \xi, \eta, \nu) =0,\qquad (\pi/2\leq\theta\leq\pi,0\leq\phi\leq 2\pi), 
		\end{equation}
		\begin{equation}\label{54b}
			\Psi^d( 0,\xi, \eta, \nu) =0,\qquad (0\leq\theta\leq\pi/2,0\leq\phi\leq 2\pi). 
		\end{equation}
	\end{subequations}
	
	Similarly from Eq.~(\ref{46}), we have
	
	\begin{equation*}
		q_1^H=[q_1^x,q_1^y]=[P(z),Q(z)]\exp{[\sigma t+i(lx+my)]},
	\end{equation*}
	
	where
	
	\begin{equation*}
		P(z)=\int_0^{4\pi}\Psi^d(z,\xi,\eta,\nu)\xi d\Omega,\quad Q(z)=\int_0^{4\pi}\Psi^d(z,\xi,\eta,\nu)\eta d\Omega.
	\end{equation*}
	
	The linear stability equations become
	
	\begin{equation}\label{54}
		\left(\sigma Sc^{-1}+k^2-\frac{d^2}{dz^2}\right)\left( \frac{d^2}{dz^2}-k^2\right)W=Rk^2N,
	\end{equation}
	\begin{equation}\label{55}
		\left(\sigma Sc^{-1}+k^2-\frac{d^2}{dz^2}\right)Z(z)=\sqrt{Ta}\frac{dW}{dz}
	\end{equation}
	\begin{equation}\label{56}
		\left(\sigma+k^2-\frac{d^2}{dz^2}\right)N+V_c\frac{d}{dz}\left(T(I_s)N+n_s\mathcal{I}_1\frac{dT(I_s)}{dG}\right)-i\frac{V_cn_sT(I_s)}{q_s}(lP+mQ)=-\frac{dn_s}{dz}W,
	\end{equation}
	
	subject to the boundary conditions
	
	\begin{equation}\label{57}
		at~~~z=0,~~~~~~W=\frac{dW}{dz}=Z(z)=\frac{dN}{dz}-V_cT(I_s)N-n_sV_c\mathcal{I}_1\frac{dT(I_s)}{dG}=0.
	\end{equation}
	
	\begin{equation}\label{58}
		at~~~z=1,~~~~~~W=\frac{d^2W}{dz^2}=\frac{dZ(z)}{dz}=\frac{dN}{dz}-V_cT(I_s)N-n_sV_c\mathcal{I}_1\frac{dT(I_s)}{dG}=0.
	\end{equation}
	
	Here, $k=\sqrt{(l^2+m^2)}$ is the non-dimensional wavenumber.
	
	Eq.~(\ref{56}) becomes (writing D = d/dz)
	
	\begin{equation}\label{58}
		\Xi_0(z)+\Xi_1(z)\int_1^z\Theta dz+(\sigma+k^2+\Xi_2(z))\Theta+\Xi_3(z)D\Theta-D^2\Theta=-Dn_sW, 
	\end{equation}
	
	where
	
	\begin{subequations}
		\begin{equation}\label{59a}
			\Xi_0(z)=V_cD\left(n_s\mathcal{I}^d\frac{dT(I_s)}{dI}\right)-i\frac{V_cn_sT(I_s)}{q_s}(lP+mQ),
		\end{equation}
		\begin{equation}\label{59b}
			\Xi_1(z)=\kappa V_cD\left(n_sI_s^c\frac{dT(I_s)}{dI}\right)
		\end{equation}
		\begin{equation}\label{59c}
			\Xi_2(z)=2\kappa V_c n_s I_s^c\frac{dT(I_s)}{dG}+V_c\frac{dT(I_s)}{dI}DI_s^d,
		\end{equation}
		\begin{equation}\label{59d}
			\Xi_3(z)=V_cT(I_s).
		\end{equation}
	\end{subequations}
	
	Introducing the new variable
	
	\begin{equation}\label{60}
		\Phi(z)=\int_1^zN(z')dz',
	\end{equation}
	
	Eqs.~\ref{54} - \ref{56} becoems
	
	\begin{equation}\label{62}
		\left(\sigma Sc^{-1}+k^2-D^2\right)\left( D^2-k^2\right)W=Rk^2D\Phi,
	\end{equation}
	
	\begin{equation}\label{63}
		\left(\sigma Sc^{-1}+k^2-D^2\right)Z(z)=\sqrt{Ta}DW
	\end{equation}
	
	\begin{equation}\label{64}
		\Xi_0(z)+\Xi_1(z)\Phi+(\sigma+k^2+\Xi_2(z))D\Phi+\Xi_3(z)D^2\Phi-D^3\Phi=-Dn_sW, 
	\end{equation}
	
	with the boundary conditions,
	
	\begin{equation}\label{65}
		at~~~z=0,~~~~~~W=DW=Z(z)=D^2\Phi-V_cT(I_s)\Phi-n_sV_c\mathcal{I}_1\frac{dT(I_s)}{dG}=0.
	\end{equation}
	
	\begin{equation}\label{66}
		at~~~z=1,~~~~~~W=D^2W=DZ(z)=D^2\Phi-V_cT(I_s)\Phi-n_sV_c\mathcal{I}_1\frac{dT(I_s)}{dG}=0,
	\end{equation}
	
	and the additional boundary condition is
	
	\begin{equation}\label{67}
		at~~~z=1,~~~~~~\Phi(z)=0.
	\end{equation}
	
	\section{SOLUTION PROCEDURE}
	
	To solve Equations (\ref{62}) and (\ref{64}), we employ a numerical approach that utilizes a fourth-order accurate finite-difference scheme based on NRK (Newton Raphson Kantorovich) iterations~\cite{19cash1980}. This numerical method enables the computation of the growth rate Re$(\sigma)$ and the generation of neutral stability curves in the $(k, R)$-plane. The neutral stability curve, denoted as $R^{(n)}(k)$ (where $n = 1, 2, 3, \ldots$), comprises an infinite number of branches. It signifies a unique solution to the linear stability problem corresponding to specific fixed parameter values. Among these branches, the most significant one is the branch where $R$ reaches its minimum value. The bioconvective solution associated with this branch, i.e., $(k_c, R_c)$, is termed the most unstable solution. This solution represents a key aspect of the bioconvection behavior and provides insights into the critical conditions for instability.
	
	
	\section{NUMERICAL RESULTS}
	
	In this article, we have conducted a systematic exploration of the impact of rotation under fixed angles of incidence. Throughout this analysis, the parameters $Sc$ and $G_t$ have been maintained as constants. Our focus has been on observing how specific constant parameter values influence the onset of bioconvection. Two constants, $Sc=20$ and $G_t=1$, have been consistently used in this study. We have concentrated on a discrete set of fixed parameter values that affect the onset of bioconvection. The parameters considered include the scattering albedo $\omega$ in the range of [0, 1], the extinction coefficient $\kappa$ with values of 0.5 and 1.0, and the cell swimming speed $V_c$ with options of 10, 15, and 20. It's noteworthy that the scattering albedo $\omega$ has been deliberately chosen to ensure that the maximum cell concentration in the basic state occurs around the midpoint of the domain's height ($z\approx 1/2$) for $\theta_i=0$. Throughout our analysis, we have consistently employed the critical intensity value $I_c=1$.
	
	
	\subsection{WHEN SCATTERING IS DOMINATED BY ABSORPTION}
	
	\begin{figure*}[!htbp]
		\includegraphics{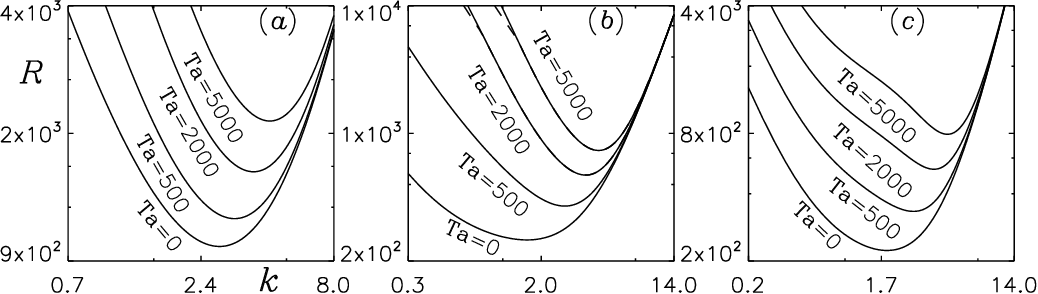}
		\caption{\label{fig3}The neutral stability curves for (a) $\theta_i=0$, (b) $\theta_i=40$, (c) $\theta_i=80$. Here, the other governing parameter values are $Sc=20, V_c=15,k=0.5$, and $\omega=0.475$. The solid line indicates the stationary branch and the dashed line indicates the oscillatory branch of the neutral cure.}
	\end{figure*}
	
	(\romannumeral 1) When extinction coefficient $\kappa=0.5$
	
	Figure~\ref{fig3} provides insight into the neutral stability curves for the Taylor number $Ta$, progressing in increasing order up to 5000, across three different scenarios ($\theta_i=0,40,80$). These simulations maintain fixed parameters of $V_c=15$, $\kappa=0.5$, and $\omega=0.475$.
	
	In the first scenario ($\theta_i=0$), as the Taylor number escalates, the critical Rayleigh number rises while the critical wavelength diminishes. Moreover, the linear stability analysis consistently reveals a stationary solution, implying that the perturbation to the basic state remains stationary in this case.
	
	Moving on to the second case ($\theta_i=40$), when the Taylor number $Ta$ spans from 0 to 1000, the most critical Rayleigh number aligns with the stationary branch of the neutral curve. This indicates the prevalence of a stationary perturbation to the basic state. Upon reaching $Ta=2000$, an oscillatory branch bifurcates from the stationary branch at $k=0.89$, persisting throughout $k\leq 0.89$. Notably, despite this oscillatory branch split, the most unstable solution continues to reside on the stationary branch of the neutral curve, affirming a stationary perturbation to the basic state. The neutral curve's characteristics, here, mirror those of an oscillatory branch emerging but the minimum Rayleigh number still appears on the stationary branch. Consequently, the perturbation to the basic state remains stationary throughout, and the critical Rayleigh number (and wavelength) demonstrate an increase as the Taylor number grows.
	
	In the third scenario ($\theta_i=80$), behaviors akin to the first case ($\theta_i=0$) are witnessed. Specifically, the critical Rayleigh number ascends while the wavelength diminishes with increasing Taylor number $Ta$ from 0 to 5000. Furthermore, a stationary perturbation is consistently observed.

	\begin{table}[!htbp]
		\caption{\label{tab1}The quantitative values of bio-convective solutions for $V_c=15$. Here for $\kappa=0.5$, $\omega=0.475$ and for $\kappa=1$, $\omega=0.61$ are fixed. The single dagger dagger symbol indicates the existence of an oscillatory branch of the neutral cure.}
		\begin{ruledtabular}
			\begin{tabular}{cccccccc}
				\multirow{2}{*}{$\theta_i$}&\multirow{2}{*}{$T_a$}& \multicolumn{3}{c}{$\kappa=0.5$} & \multicolumn{3}{c}{$\kappa=1$}\\\cline{3-5}\cline{6-8} & & $\lambda_c$ & $R_c$ & $Im(\sigma)$ &$ \lambda_c$ & $R_c$ & $Im(\sigma)$\\
				\hline
				
				0 & 0  & 2.23 & 979.51 & 0 & 2.10 & 559.72 & 0 \\
				0 & 100 & 2.14 & 1019.07 & 0 & 2.02 & 582.88 & 0 \\
				0 & 500 & 1.95 & 1152.55 & 0 & 1.86 & 660.40 & 0 \\
				0 & 1000 & 1.81 & 1289.58 & 0 & 1.71 & 739.20 & 0 \\
				0 & 2000 & 1.64 & 1516.08 & 0 & 1.55 & 868.84 & 0 \\
				0 & 5000 & 1.41 & 2039.69 & 0 & 1.34$^\dagger$ & 1167.13 & 0 \\
				
				40 & 0  & 3.83 & 218.03 & 0 & 2.13$^\dagger$ & 276.87 & 0 \\
				40 & 100 & 3.06 & 266.39 & 0 & 2.04$^\dagger$ & 301.02 & 0 \\
				40 & 500 & 2.18 & 397.97 & 0 & 1.77$^\dagger$ & 379.70 & 0 \\
				40 & 1000 & 1.87 & 514.08 & 0 & 1.61$^\dagger$ & 457.86 & 0 \\
				40 & 2000 & 1.61$^\dagger$ & 690.53 & 0 & 1.43$^\dagger$ & 583.30 & 0 \\
				40 & 5000 & 1.33$^\dagger$ & 1069.92 & 0 & 1.21$^\dagger$ & 861.48 & 0 \\
				
				80 & 0  & 3.44 & 182.60 & 0 & 2.57 & 267.94 & 0 \\
				80 & 100 & 2.98 & 211.85 & 0 & 2.35 & 294.01 & 0 \\
				80 & 500 & 2.25 & 298.96 & 0 & 1.94 & 372.10 & 0 \\
				80 & 1000 & 1.91 & 380.60 & 0 & 1.73 & 444.43 & 0 \\
				80 & 2000 & 1.63 & 508.75 & 0 & 1.52 & 555.89 & 0 \\
				80 & 5000 & 1.30 & 791.08 & 0 & 1.24 & 795.54 & 0 \\
			\end{tabular}    
		\end{ruledtabular}
	\end{table}	   
	
	
	(\romannumeral 2) When extinction coefficient $\kappa=1$
	
	\begin{figure*}[!htbp]
		\includegraphics{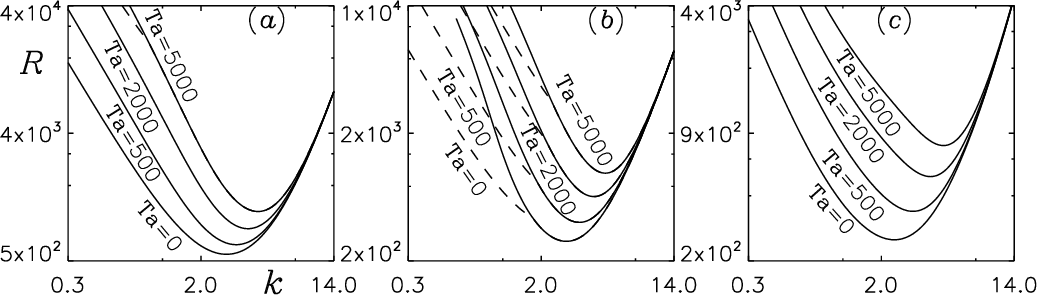}
		\caption{\label{fig4}The neutral stability curves for (a) $\theta_i=0$, (b) $\theta_i=40$, (c) $\theta_i=80$. Here, the other governing parameter values are $Sc=20, V_c=15,k=1$, and $\omega=0.61$. The solid line indicates the stationary branch and the dashed line indicates the oscillatory branch of the neutral cure.}
	\end{figure*}

	In Figure~\ref{fig4}(a), the neutral stability curves for the Taylor number $Ta$ are depicted in ascending order up to 5000, with $\theta_i=0$ and fixed parameters $V_c=15$, $\kappa=1$, and $\omega=0.61$. At $Ta=0$, the linear stability analysis indicates a stationary solution. As the Taylor number progresses up to 2000, the solution maintains its stationary nature while the critical Rayleigh number increases and the critical wavelength decreases. When $Ta$ reaches 5000, an oscillatory branch bifurcates from the stationary branch; however, the stationary branch retains the minimum Rayleigh number.
	
	In Figure~\ref{fig4}(b), the neutral curves are presented for $\theta_i=40$. For all Taylor number values ranging from 0 to 5000, an oscillatory branch emerges from the neutral curve's stationary branch. Nevertheless, the most unstable solution remains situated on the stationary branch, resulting in a stationary bifurcation from the basic state. With increasing Taylor numbers, the critical Rayleigh number increases while the wavelength decreases.
	
	Figure~\ref{fig4}(c) portrays the neutral curves for $\theta_i=80$. In this scenario, the linear stability analysis predicts solely stationary perturbations to the basic state for all Taylor number values. Here, an increase in the Taylor number corresponds to an increase in the critical Rayleigh number and a decrease in wavelength.
	
	The critical Rayleigh numbers and wavelengths corresponding to $V_c=15$ are summarized in Table~\ref{tab1}.
	
	
	Here also shows the numerical observations for $V_c=10$ and $V_c=20$ in the Table~\ref{tab2} and Table~\ref{tab3}.

	\begin{table}[h]
		\caption{\label{tab2}The quantitative values of bio-convective solutions for $V_c=10$. Here for $\kappa=0.5$, $\omega=0.475$ and for $\kappa=1$, $\omega=0.6$ are fixed. The single dagger symbol indicates the existence of an oscillatory branch of the neutral cure.}
		\begin{ruledtabular}
			\begin{tabular}{cccccccc}
				\multirow{2}{*}{$\theta_i$}&\multirow{2}{*}{$T_a$}& \multicolumn{3}{c}{$\kappa=0.5$} & \multicolumn{3}{c}{$\kappa=1$}\\\cline{3-5}\cline{6-8} & & $\lambda_c$ & $R_c$ & $Im(\sigma)$ &$ \lambda_c$ & $R_c$ & $Im(\sigma)$\\
				\hline
				
				0 & 0  & 2.14 & 1985.17 & 0 & 2.27 & 778.29 & 0 \\
				0 & 100 & 2.09 & 2042.02 & 0 & 2.14 & 814.20 & 0 \\
				0 & 500 & 1.94 & 2244.66 & 0 & 1.92 & 932.76 & 0 \\
				0 & 1000 & 1.80 & 2462.30 & 0 & 1.80 & 1052.06 & 0 \\
				0 & 2000 & 1.65 & 2832.57 & 0 & 1.61 & 1247.15 & 0 \\
				0 & 5000 & 1.44 & 3705.79 & 0 & 1.39 & 1695.11 & 0 \\
				
				40 & 0  & 5.26 & 390.95 & 0 & 3.15 & 226.72 & 0 \\
				40 & 100 & 3.47 & 493.50 & 0 & 2.69 & 266.05 & 0 \\
				40 & 500 & 2.32 & 721.25 & 0 & 2.09 & 380.76 & 0 \\
				40 & 1000 & 1.99 & 906.29 & 0 & 1.83$^\dagger$ & 486.58 & 0 \\
				40 & 2000 & 1.71 & 1180.37 & 0 & 1.57$^\dagger$ & 650.52 & 0 \\
				40 & 5000 & 1.43 & 1765.35 & 0 & 1.32$^\dagger$ & 1007.21 & 0 \\
				
				80 & 0  & 5.22 & 148.21 & 0 & 3.44 & 177.88 & 0 \\
				80 & 100 & 3.93 & 187.03 & 0 & 2.98 & 207.12 & 0 \\
				80 & 500 & 2.57 & 302.88 & 0 & 2.25 & 294.87 & 0 \\
				80 & 1000 & 2.09 & 411.07 & 0 & 1.94 & 377.84 & 0 \\
				80 & 2000 & 1.73 & 580.94 & 0 & 1.65 & 509.06 & 0 \\
				80 & 5000 & 1.37 & 955.68 & 0 & 1.33 & 800.16 & 0 \\
			\end{tabular}    
		\end{ruledtabular}
	\end{table}
	
	\begin{table}[h]
		\caption{\label{tab3}The quantitative values of bio-convective solutions for $V_c=20$. Here for $\kappa=0.5$, $\omega=0.48$ and for $\kappa=1$, $\omega=0.612$ are fixed. The single dagger symbol indicates the existence of an oscillatory branch of the neutral cure and the double dagger symbol indicates an overstable solution.}
		\begin{ruledtabular}
			\begin{tabular}{cccccccc}
				\multirow{2}{*}{$\theta_i$}&\multirow{2}{*}{$T_a$}& \multicolumn{3}{c}{$\kappa=0.5$} & \multicolumn{3}{c}{$\kappa=1$}\\\cline{3-5}\cline{6-8} & & $\lambda_c$ & $R_c$ & $Im(\sigma)$ &$ \lambda_c$ & $R_c$ & $Im(\sigma)$\\
				\hline
				
				0 & 0  & 2.27 & 658.84 & 0 & 2.10 & 392.47 & 0 \\
				0 & 100 & 2.14 & 689.13 & 0 & 1.99 & 412.43 & 0 \\
				0 & 500 & 1.92 & 789.34 & 0 & 1.77 & 477.06 & 0 \\
				0 & 1000 & 1.77 & 890.37 & 0 & 1.64 & 540.80 & 0 \\
				0 & 2000 & 1.61 & 1055.63 & 0 & 1.48 & 643.63 & 0 \\
				0 & 5000 & 1.37$^\dagger$ & 1435.01 & 0 & 1.29 & 875.89 & 0 \\
				
				40 & 0  & 2.69 & 225.15 & 0 & 2.34$^\ddagger$ & 409.74$^\ddagger$ & 9.34 \\
				40 & 100 & 2.45 & 257.38 & 0 & 1.61$^\dagger$ & 430.73 & 0 \\
				40 & 500 & 1.95$^\dagger$ & 354.77 & 0 & 1.48$^\dagger$ & 491.03 & 0 \\
				40 & 1000 & 1.71$^\dagger$ & 446.39 & 0 & 1.40$^\dagger$ & 554.09 & 0 \\
				40 & 2000 & 1.50$^\dagger$ & 589.34 & 0 & 1.29$^\dagger$ & 658.71 & 0 \\
				40 & 5000 & 1.26$^\dagger$ & 900.03 & 0 & 1.12$^\dagger$ & 896.58 & 0 \\
				
				80 & 0  & 2.35$^\dagger$ & 247.63 & 0 & 2.05 & 428.37 & 0 \\
				80 & 100 & 2.17$^\dagger$ & 274.60 & 0 & 1.94 & 453.30 & 0 \\
				80 & 500 & 1.83$^\dagger$ & 359.39 & 0 & 1.68 & 528.31 & 0 \\
				80 & 1000 & 1.64$^\dagger$ & 441.21 & 0 & 1.52 & 597.19 & 0 \\
				80 & 2000 & 1.46$^\dagger$ & 570.28 & 0 & 1.37 & 701.65 & 0 \\
				80 & 5000 & 1.23$^\dagger$ & 852.29 & 0 & 1.15 & 920.83 & 0 \\
			\end{tabular}    
		\end{ruledtabular}
	\end{table}
	
	
	\begin{figure*}[!htbp]
		\includegraphics{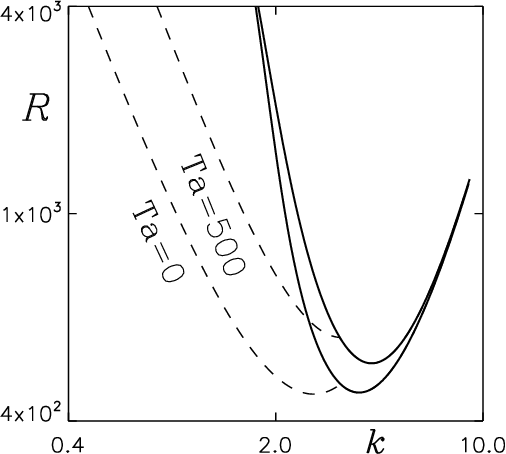}
		\caption{\label{fig5}The existence of overstable solution at $V_c=20$, $\kappa=1$, $\omega=0.612$, and $\theta_i=40$. The solid line indicates the stationary branch and the dashed line indicates the oscillatory branch of the neutral cure.}
	\end{figure*}
	
	\subsubsection{$THE$ $OVERSTABLE$ $SOLUTION$}

	Overstable solutions have been identified for specific parameter configurations. In cases where overstable solutions occur, a singular oscillatory segment of the neutral curve intersects the stationary segment at a critical wavenumber $k = k_b$, persisting for wavenumbers $k \leq k_b$. Consequently, an oscillatory segment emerges from the fundamental state at the critical wavenumber $k = k_b$. The numerical linear stability analysis has been conducted for various parameter values, specifically $ V_ =$ 10, 15, 20, and $\kappa =$ 0.5, 1.0.
	
	No overstable solutions were discovered for $\kappa =$ 0.5 and $V_c =$ 10, 15, 20, nor for $\kappa =$ 1.0 and $V_c =$ 10, 15.
	For the case of $V_c =$ 20, $\kappa =$ 1.0, and $\omega =$ 0.612, when the Taylor number $Ta$ is set to 0, an oscillatory segment of the neutral curve intersects the stationary segment at $k = 3.31$, persisting for wavenumbers $k \leq 3.31$. However, on this stationary segment, the most unstable solution remains $(R_c,\lambda_c) = (409.74,2.69)$, with a real growth rate $Re(\sigma) = 9.34$ (as shown in Fig.~\ref{fig5}). As the Taylor number $Ta$ is increased to 100 or 500, an oscillatory branch becomes linked to the stationary branch around $k = 3.35$. Nevertheless, despite this connection, the most unstable mode is situated on the stationary branch, leading to a stationary solution.
	
	
	\section{Conclusion}
	
	This article focuses on investigating the influence of rotation on the initiation of phototactic bioconvection within a suspension of isotropic scattering phototactic algae. The suspension is subjected to oblique collimated irradiation from above. The study employs linear perturbation theory to analyze the linear stability of this suspension.
	
	Under constant parameter conditions, the linear stability analysis of the steady state reveals a transition from stationary to oscillatory solutions and vice versa as the Taylor number or rotation rate increases. The emergence of an overstable solution is attributed to the conflict between stabilizing and destabilizing processes. As the Taylor number rises, the critical Rayleigh number also increases, accompanied by a reduction in pattern wavelength.
	
	Nevertheless, it's important to acknowledge that the proposed phototaxis model should ideally be validated through quantitative experimental data on bioconvection in a purely phototactic algal suspension. Presently, such data are unavailable. Therefore, efforts should be directed toward identifying suitable microorganism species predominantly exhibiting phototaxis. Currently, naturally occurring algal species tend to possess either gravitactic or gyrotactic traits alongside phototaxis. It's worth noting that the proposed model could also be applied to simulate phototactic bioconvection within an anisotropic algal suspension, offering insights into various intriguing problems.  
	
	
	\section*{Data Availability}
	The data that support the plots within this paper and other findings of this study are available in the article.

	\nocite{*}
	\bibliography{aipsamp}
	
\end{document}